\documentclass[a4paper]{article}
\usepackage{amssymb}
\usepackage{amsfonts}
\usepackage{amsmath}

\input{tcilatex}

\begin{document}

\title{The multi-degree of coverings on Lie groups}
\author{Haibao Duan\thanks{%
Supported by NSFC 11131008; 11661131004} \\
Institute of Mathematics, Chinese Academy of Sciences, \\
Beijing 100190, P.R. China \and Siye Wu \\
Department of Mathematics, National Tsing Hua University, \\
Hsinchu 30013, Taiwan}
\maketitle

\begin{abstract}
We associate to each covering map of simple Lie groups a sequence of
integers, called \textsl{the multi--degree of the covering}; develop a
unified method to evaluate the sequence; and apply the results to solve two
extension problems of mathematical physics.

\begin{description}
\item[2010 Mathematical Subject Classification: ] 57T10; 55R20;  81T13.

\item[Key words and phrases:] cohomology of Lie groups; Serre spectral
sequence; Schubert calculus; WZW models; Gauge theory

\item[Emails:] dhb@math.ac.cn; swu@math.nthu.edu.tw
\end{description}
\end{abstract}

\section{Main results}

Covering maps on Lie groups are essential to extend the constructions and
calculations of various physical models associated with simply connected Lie
groups, to that associated with non--simply connected Lie groups. Naturally,
topological invariants of the coverings are useful to formulate, or to
solve, the relevant extension problems. In this paper we introduce for each
covering map of Lie groups an invariant, called \textsl{the multi--degree of
the covering}; extend Schubert calculus to evaluate the invariant; and apply
the results to two outstanding topological problems arising from the studies
of the Wess--Zumino--Witten models and the topological Gauge theories \cite%
{DE,FGS,MR}. The main tool in our approach\textsl{\ }is \textsl{the Chow
rings} of Lie groups, introduced by Grothendieck \cite{G} in 1958.

To clarify the concept of \textsl{multi--degree }of a covering on Lie groups
we recall the classical results of Hopf and Chevalley on the real cohomology
of Lie groups. The Lie groups $G$ under consideration will be compact and
connected. For a maximal torus $T$ on $G$ the dimension $n:=\dim T$ is an
invariant of $G$, called \textsl{the rank of }$G$. The multiplication $\mu
:G\times G\rightarrow G$ on $G$ defines the \textsl{co--product }on the real
cohomology algebra $H^{\ast }(G;\mathbb{R})$

\begin{enumerate}
\item[(1.1)] $\mu ^{\ast }:H^{\ast }(G;\mathbb{R})\rightarrow H^{\ast
}(G\times G;\mathbb{R})\cong H^{\ast }(G;\mathbb{R})\otimes H^{\ast }(G;%
\mathbb{R})$,
\end{enumerate}

\noindent where the isomorphism is given by the Kunneth formula. An element $%
y\in H^{\ast }(G;\mathbb{R})$ is called \textsl{primitive} if the relation $%
\mu ^{\ast }(y)=y\otimes 1+1\otimes y$ is satisfied. Since the linear
combinations of primitive elements are also primitive such elements form a
subspace $P(G;\mathbb{R})$ of the algebra $H^{\ast }(G;\mathbb{R})$. Hopf 
\cite{H} has shown that

\bigskip

\noindent \textbf{Theorem 1.1.} \textsl{For a Lie group }$G$\textsl{\ with
rank }$n$\textsl{\ we have }$\dim P(G;\mathbb{R})=n$\textsl{. }

\textsl{Moreover, if }$\left\{ y_{1},\cdots ,y_{n}\right\} $\textsl{\ is a
basis of the space }$P(G;\mathbb{R})$\textsl{, then the real cohomology
algebra} $H^{\ast }(G;\mathbb{R})$ \textsl{is the exterior algebra generated
by} $y_{1},\cdots ,y_{n}$\textsl{:}

\begin{enumerate}
\item[(1.2)] $H^{\ast }(G;\mathbb{R})=\Lambda _{\mathbb{R}}(y_{1},\cdots
,y_{n})$\textsl{,} $\deg y_{i}\equiv 1\func{mod}2$.$\square $
\end{enumerate}

By (1.2) the sequence $I_{G}=\{r_{1},\cdots ,r_{n}\}$ of integers defined by

\begin{quote}
$r_{i}:=\frac{1}{2}(\deg y_{i}+1),$ $1\leq i\leq n,$
\end{quote}

\noindent is also an invariant of the group $G$, which has been shown by
Chevalley \cite{Ch0},\cite[(3.2)Theorem]{R} to be the degree sequence of 
\textsl{the basic Weyl invariants of }$G$. We may assume that the sequence $%
I_{G}$ is ordered by $r_{1}\leq \cdots \leq r_{n}$. In particular, if the
group $G$ is simple one has $2=r_{1}<\cdots <r_{n}$ (see in Table 1 ).

In general, for a CW--complex $X$ the torsion subgroup $\tau (X)$ of the
integral cohomology $H^{\ast }(X)$ is an ideal and therefore, defines the
(graded) quotient ring $\mathcal{F}(X):=H^{\ast }(X)/\tau (X)$. In view of
the obvious additive decomposition

\begin{quote}
$H^{\ast }(X)=\mathcal{F}(X)\oplus \tau (X)$
\end{quote}

\noindent we may call $\mathcal{F}(X)$ \textsl{the torsion free part }of the
integral cohomology\textsl{\ }$H^{\ast }(X)$. Furthermore, if $%
f:X\rightarrow Y$ is a continuous map between two CW--complexes, the induce
ring map $f^{\ast }$ on the cohomologies preserves the torsion ideals.
Therefore, passing to the quotients yields the (graded) ring map $f^{\#}:%
\mathcal{F}(Y)\rightarrow \mathcal{F}(X)$. In particular, the multiplication 
$\mu $ on $G$ furnishes the ring $\mathcal{F}(G)$ with a co--product

\begin{quote}
$\mu ^{\#}:\mathcal{F}(G)\rightarrow \mathcal{F}(G)\otimes \mathcal{F}(G)$.
\end{quote}

\noindent The following result may be seen as an integral refinement of
Theorem 1.1.

\bigskip

\noindent \textbf{Theorem A.} \textsl{For a Lie group }$G$\textsl{\ with }$%
I_{G}=\{r_{1},\cdots ,r_{n}\}$\textsl{, there exist }$n$ \textsl{elements }$%
x_{1},\cdots ,x_{n}\in \mathcal{F}(G)$\textsl{,} $\deg x_{i}=2r_{i}-1$%
\textsl{,\ so that}

\begin{quote}
\textsl{i) }$\mathcal{F}(G)=\Lambda (x_{1},\cdots ,x_{n})$\textsl{\ (i.e.
the exterior ring over }$\mathbb{Z}$\textsl{)}

\textsl{ii) }$\mu ^{\#}(x_{i})=x_{i}\otimes 1+1\otimes x_{i}$\textsl{, }$%
1\leq i\leq n$\textsl{.}
\end{quote}

\noindent \textsl{In addition,, if the group }$G$\textsl{\ is simple,} 
\textsl{then the generators }$x_{1},\cdots ,x_{n}$ \textsl{satisfying i) and
ii)} \textsl{are} \textsl{unique up to the sign }$\pm $\textsl{\ convention.}

\bigskip

Let $G$ be a simple Lie group with non--trivial center $\mathcal{Z}(G)$. The
quotient $c:G\rightarrow PG:=G/\mathcal{Z}(G)$ is a covering map of Lie
groups, and induces a ring map on the torsion free parts of the cohomologies

\begin{enumerate}
\item[(1.3)] $c^{\#}:\mathcal{F}(PG)=\Lambda (y_{1},\cdots
,y_{n})\rightarrow \mathcal{F}(G)=\Lambda (x_{1},\cdots ,x_{n})$,
\end{enumerate}

\noindent where $\deg y_{i}=\deg x_{i}$ by the well known fact $I_{G}=I_{PG}$
of invariant theory. Granted with Theorem A we shall show that

\bigskip

\noindent \textbf{Theorem B. }\textsl{There exists a unique sequence }$%
\{a_{1},\cdots ,a_{n}\}$\textsl{\ of }$n$ \textsl{positive integers, so that}

\textsl{i) the map }$c^{\#}$\textsl{\ is given by }$c^{\#}(y_{i})=a_{i}\cdot
x_{i}$\textsl{, }$1\leq i\leq n$\textsl{;}

\textsl{ii) the product }$a_{1}\cdot \cdots \cdot a_{n}$\textsl{\ is the
order }$\left\vert \mathcal{Z}(G)\right\vert $\textsl{\ of the center.}

\bigskip

The unique sequence $\{a_{1},\cdots ,a_{n}\}$ obtained by Theorem B, written 
$\mathcal{D}(G,PG)$ and called \textsl{the multi--degree} of the covering $c$%
, will be the main concern of this work. To be precise we tabulate all the
simply connected simple Lie groups $G$ with non--trivial centers $\mathcal{Z}%
(G)$, together with the degree sequences $I_{G}$ of their basic Weyl
invariants, using the table below.

\begin{center}
{\normalsize Table 1. The simply connected simple Lie groups with
non--trivial centers}

\begin{tabular}{|l|l|l|l|}
\hline
${\small G}$ & rank$G$ & $\qquad I_{G}$ & $\mathcal{Z}(G)$ \\ \hline
${\small SU(n)}$ & $n-1$ & $\left\{ 2,3,\cdots ,n\right\} $ & $\mathbb{Z}%
_{n} $ \\ \hline
${\small Sp(n)}$ & $n$ & $\left\{ 2,4,\cdots ,2n\right\} $ & $\mathbb{Z}_{2}$
\\ \hline
${\small Spin(2n+1)}$ & $n$ & $\left\{ 2,4,\cdots ,2n\right\} $ & $\mathbb{Z}%
_{2}$ \\ \hline
${\small Spin(2n)},n\equiv 1\func{mod}2$ & $n$ & $\left\{ 2,4,\cdots
,2(n-1)\}\amalg \{n\right\} $ & $\mathbb{Z}_{4}$ \\ \hline
${\small Spin(2n)},n\equiv 0\func{mod}2$ & $n$ & $\left\{ 2,4,\cdots
,2(n-1)\right\} \amalg \{n\}$ & $\mathbb{Z}_{2}{\small \oplus }\mathbb{Z}%
_{2} $ \\ \hline
${\small E}_{6}$ & $6$ & $\left\{ 2,5,6,8,9,12\right\} $ & $\mathbb{Z}_{3}$
\\ \hline
${\small E}_{7}$ & $7$ & $\left\{ 2,6,8,10,12,14,18\right\} $ & $\mathbb{Z}%
_{2}$ \\ \hline
\end{tabular}
\end{center}

Let $J(PG)$ be the subring of $H^{\ast }(PG)$\ generated multiplicatively by
the second cohomology group $H^{2}(PG)$, together with the multiplicative
unit $1\in H^{0}(PG)$. By computing with the Chow rings \cite{G} of the
adjoint Lie groups $PG$ we shall show that

\bigskip

\noindent \textbf{Theorem C.} \textsl{For each simply connected Lie group }$%
G $ \textsl{in Table 1 we have}

\begin{enumerate}
\item[(1.4)] $H^{2}(PG)=\mathcal{Z}(G)$.
\end{enumerate}

\noindent \textsl{Moreover, there exist generators}

\begin{quote}
$\omega \in H^{2}(PG)$ \textsl{for} $G\neq Spin(2n)$\textsl{\ with }$n\equiv
0\func{mod}2$\textsl{;\ or}

$\omega _{1}$\textsl{, }$\omega _{2}\in H^{2}(PG)=\mathbb{Z}_{2}\oplus 
\mathbb{Z}_{2}$ \textsl{for} $G=Spin(2n)$\textsl{\ with }$n\equiv 0\func{mod}%
2$
\end{quote}

\noindent \textsl{so that the ring }$J(PG)$\textsl{\ has the following
presentation:}

\begin{quote}
\textsl{i)} $J(PSU(n))=\frac{\mathbb{Z}[\omega ]}{\left\langle b_{n,r}\omega
^{r},\text{ }1\leq r\leq n\right\rangle },$\ \textsl{where} $b_{n,r}=g.c.d.\{%
\binom{n}{1},\cdots ,\binom{n}{r}\}$;

\textsl{ii)} $J(PSp(n))=\frac{\mathbb{Z}[\omega ]}{\left\langle 2\omega
,\omega ^{2^{r+1}}\right\rangle },$ \textsl{where} $n=2^{r}(2k+1)$;

\textsl{iii)} $J(PSpin(2n+1))=\frac{\mathbb{Z}[\omega ]}{\left\langle
2\omega ,\omega ^{2^{s+1}}\right\rangle }$\textsl{,} \textsl{where} $%
2^{s}\leq n<2^{s+1}$;

\textsl{iv)} $J(PSpin(2n))=\left\{ 
\begin{array}{c}
\frac{\mathbb{Z}[\omega ]}{\left\langle 4\omega ,2\omega ^{2},\omega
^{2^{r+1}}\right\rangle }\text{ \textsl{if} }n\equiv 1\func{mod}2\text{; }
\\ 
\frac{\mathbb{Z}[\omega _{1},\omega _{2}]}{\left\langle 2\omega _{1},2\omega
_{2},\omega _{1}^{2^{t}},\omega _{2}^{2^{r+1}}\right\rangle }\text{ \textsl{%
if} }n\equiv 0\func{mod}2\text{,}%
\end{array}%
\right. $

\textsl{where} $n=2^{t}(2k+1)${\small , }$2^{r}<n\leq 2^{r+1}$;

\textsl{v)} $J(PE_{6})=\frac{\mathbb{Z}[\omega ]}{\left\langle 3\omega
,\omega ^{9}\right\rangle }$; \quad \textsl{vi)} $J(PE_{7})=\frac{\mathbb{Z}%
[\omega ]}{\left\langle 2\omega ,\omega ^{2}\right\rangle }$.
\end{quote}

Our main result implies that, the multi-degree sequence $\mathcal{D}(G,PG)$
is determined entirely by the structure of the ring $J(PG)$ given by Theorem
C (see Theorem 4.4).

\bigskip

\noindent \textbf{Theorem D. }\textsl{For a Lie group }$G$\textsl{\ in Table
1} \textsl{the multi--degree }$\mathcal{D}(G,PG)$\textsl{\ is}

\begin{center}
\begin{tabular}{|l|l|}
\hline
${\small G}$ & $\mathcal{D}(G,PG)$ \\ \hline\hline
${\small SU(n)}$ & ${\small \{a}_{1}{\small ,}\cdots {\small ,a}_{n-1}%
{\small \},}$\textsl{\ }${\small a}_{k}{\small =}\frac{b_{n,k}}{b_{n,k+1}}%
{\small ,}$ \\ \hline
${\small Sp(n)},$ ${\small n=2}^{t}{\small (2b+1)}$ & ${\small \{1,\cdots
,1,2}_{(2^{t})}{\small ,1\cdots ,1\}}$ \\ \hline
${\small Spin(2n+1)}${\small , }${\small 2}^{s}{\small \leq n<2}^{s+1}$ & $%
{\small \{1,\cdots ,1,2}_{(2^{s})}{\small ,1\cdots ,1\}}$ \\ \hline
${\small Spin(2n)},{\small 2}^{s}{\small <n=2b+1\leq 2}^{s+1}$ & ${\small %
\{2,1,\cdots ,1,2}_{(2^{s})}{\small ,1,\cdots ,1\}}$ \\ \hline
${\small Spin(2n)},{\small 2}^{s}{\small <n=2}^{t}{\small (2b+1)\leq 2}^{s+1}%
{\small ,t\geq 1}$ & ${\small \{1,\cdots ,1,2}_{(2^{t-1})}{\small ,1,\cdots
,1,2}_{(2^{s})}{\small ,1,\cdots ,1\}}$ \\ \hline
${\small E}_{6}$ or ${\small E}_{7}$ & ${\small \{1,1,1,1,3,1\}}$ or $%
{\small \{2,1,1,1,1,1,1\}}$ \\ \hline
\end{tabular}
\end{center}

\noindent \textsl{where the notion }$2_{(h)}$\textsl{\ stands for }$a_{h}=2$%
\textsl{, and} \textsl{where if }$G=SU(n)$\textsl{\ and the integer }$n$%
\textsl{\ has the prime factorization }$n=p_{1}^{r_{1}}\cdots p_{t}^{r_{t}}$%
\textsl{, then (see \cite[Lemma 3.1]{DL})}

\begin{enumerate}
\item[(1.5)] $a_{k}=\left\{ 
\begin{tabular}{l}
$p_{i}$ \textsl{if}$\ k=p_{i}^{s}$ \textsl{with} $1\leq i\leq t$ \textsl{and}
$1\leq s\leq r_{i}$\textsl{;} \\ 
$1\text{ \textsl{otherwise}}$\textsl{.}%
\end{tabular}%
\right. $
\end{enumerate}

\bigskip

Theorem D is a natural extension of \cite[Theorem 1.1]{DL}, where the
multi-degree of the universal covering $SU(n)\rightarrow PSU(n)$ on the
projective unitary group $PSU(n)$ has been computed. On the other hand, E.B.
Dynkin \cite{Dy1,Dy2} raised the problem to determine induce actions of Lie
group homomorphisms $G\rightarrow G^{\prime }$ on the cohomologies, which is
essentially solved by Theorem D for the cases of universal coverings of the
adjoint Lie groups. Nevertheless, the present paper is motivated by two
topics from mathematical physics. To start with we note by Theorem D that

\bigskip

\noindent \textbf{Corollary 1.2.} \textsl{For a Lie group }$G$\textsl{\ in
Table 1 the leading term }$a_{1}$ \textsl{of the sequence }$\mathcal{D}%
(G,PG) $ \textsl{is either }$1$\textsl{\ or }$2$, \textsl{where }$a_{1}=2$%
\textsl{\ occurs if and only if }$G$\textsl{\ is isomorphic to one of the
following groups}

\begin{quote}
\textsl{i) }$SU(n)$ \textsl{with} $n\equiv 0\func{mod}2$\textsl{;} \quad

\textsl{ii) }$Sp(n)$ \textsl{with} $n\equiv 1\func{mod}2$\textsl{; }

\textsl{iii) }$Spin(2^{t}(2b+1))$ \textsl{with} $t=1$, $2$; \quad

\textsl{iv) }$E_{7}$.$\square $
\end{quote}

In the study on the Wess--Zumino--Witten models with simple Lie groups
Felder, Gawedzki and Kupiainen obtained the short exact sequence

\begin{enumerate}
\item[(1.6)] $0\rightarrow H_{3}(G)(=\mathbb{Z})\overset{c_{\ast }}{%
\rightarrow }H_{3}(PG)\rightarrow \mathcal{Z}(G)\rightarrow 0$,
\end{enumerate}

\noindent and investigated its extension problem \cite[Appendix 1]{FGS},
where $G$ is simply connected and simple. This problem was emphasized by
Dijkgraaf and Witten in the work \cite{DE} on the topological Gauge
theories, and by Mathai and Rosenberg \cite{MR} in the study on the
relationship between Langlands duality and T--duality for compact Lie
groups. Since the map $c_{\ast }$ in (1.6) is Kronecker dual to the map $%
c^{\ast }$ on $H^{3}(PG)$ we get from Corollary 1.2 the following result.

\bigskip

\noindent \textbf{Corollary 1.3. }\textsl{For the Lie groups }$G$\textsl{\
given in Table 1 we have}

\begin{quote}
$H_{3}(PG)=\mathbb{Z}\oplus \mathcal{Z}(G)$
\end{quote}

\noindent \textsl{with the following exceptions}

\begin{quote}
\textsl{i)} $H_{3}(PSU(2m))=$\textsl{\ }$\mathbb{Z\oplus Z}_{m}$\textsl{,
where} $m\in \mathbb{Z}$\textsl{;}

\textsl{ii)} $H_{3}(PSp(2b+1))=$\textsl{\ }$H_{3}(PE_{7})=\mathbb{Z}$\textsl{%
, where} $b\in \mathbb{Z}$\textsl{;}

\textsl{iii)} $H_{3}(PSpin(2^{t}(2b+1)))=$\textsl{\ }$\mathbb{Z\oplus Z}_{2}$%
\textsl{, where} $t=1,2,$ $b\in \mathbb{Z}\QTR{sl}{.}\square $
\end{quote}

\noindent \textbf{Remark 1.4.} In the inspiring work \cite{MR} Mathai and
Rosenberg have computed the leading terms $a_{1}$ of the sequences $\mathcal{%
D}(G,PG)$. In the cases of $G=SU(n)$ with $n$ divisible by $4$ and $%
G=Spin(2(2b+1))$, our results in Corollary 1.2 are different with the ones
stated in \cite[Theorem 1, (1), (3)]{MR}. We note that these are precisely
the cases where the integral cohomologies $H^{\ast }(G)$ have torsion
elements of the order $4$. Therefore, working with the $\mathbb{Z}_{2}$
algebra $H^{\ast }(G;\mathbb{Z}_{2})$ alone may not suffice to decide $a_{1}$%
.$\square $

\bigskip

Let $B_{G}$ be the classifying space of a simple Lie group $G$. By the
naturality of the cohomology suspension $\tau _{G}:H^{r}(B_{G})\rightarrow
H^{r-1}(G)$ \cite[p.22]{JH} in the universal $G$--bundle

\begin{quote}
$G\hookrightarrow E_{G}\rightarrow B_{G}$
\end{quote}

\noindent the covering $c:G\rightarrow PG$ induces the commutative diagram

\begin{enumerate}
\item[(1.7)] 
\begin{tabular}{lll}
$H^{4}(B_{PG})$ & $\overset{\tau _{PG}}{\rightarrow }$ & $H^{3}(PG)$ \\ 
$Bc^{\ast }\downharpoonright $ &  & $c^{\ast }\downharpoonright $ \\ 
$H^{4}(B_{G})$ & $\overset{\cong \tau _{G}}{\rightarrow }$ & $H^{3}(G)$%
\end{tabular}
(compare with \cite[(4.14)]{DE})
\end{enumerate}

\noindent where, as being pointed out in \cite{DE} that, the four vertices
groups in (1.7) are all isomorphic to $\mathbb{Z}$, and where the suspension 
$\tau _{G}$ at the bottom is an isomorphism when $G$ is simply connected.
Let $\omega $ and $\widetilde{\omega }$ denote the generators of
respectively $H^{4}(B_{PG})$ and $H^{4}(B_{G})$, $\xi $ and $\widetilde{\xi }
$ denote the generators of respectively $H^{3}(PG)$ and $H^{3}(G)$. In \cite[%
Section 4.3]{DE} Dijkgraaf and Witten raised the interesting problem to
decide the pair $(\alpha ,\beta )$ of integers characterized by the relations

\begin{quote}
$Bc^{\ast }(\omega )=\alpha \cdot \widetilde{\omega }$, $\tau _{PG}(\omega
)=\beta \cdot \xi $,
\end{quote}

\noindent where the pair $(\alpha ,\beta )$ has shown to be $(\frac{2n}{%
g.c.d.\{2,n-1\}},n)$ for the case $G=SU(n)$.

On the other hand, by the commutativity of the diagram (1.7) one has $\alpha
=a_{1}\cdot \beta $, where $a_{1}$ is the leading term of the sequence $%
\mathcal{D}(G,PG)$. Therefore, combining results of Corollary 1.2 with
Dijkgraaf--Witten's formula \cite[(4.18)]{DE} evaluating the integer $\alpha 
$ we obtain the following results.

\bigskip

\noindent \textbf{Corollary 1.5. }\textsl{For a Lie group }$G$\textsl{\ in
Table 1 the pair }$(\alpha ,\beta )$ \textsl{of integers is given by the
table:}

\begin{center}
\begin{tabular}{|l|l|}
\hline
$G$ & $(\alpha ,\beta )$ \\ \hline
$SU(n)$ & $(\frac{2n}{g.c.d.\{2,n-1\}},n)$ \\ \hline
$Spin(2n+1)$ & $(2,2)$ \\ \hline
$Sp(n),$ $n=2^{t}(2b+1)$ & $(8,4)$, $(4,4)$, $(2,2)$ or $(1,1)$ for $t=0,1,2$
or $\geq 3$ \\ \hline
$Spin(2n),$ $n=2^{t}(2b+1)$ & $(8,4)$, $(4,4)$ or $(2,2)$ for $t=0,1$ or $%
\geq 2$ \\ \hline
$E_{6}$ & $(3,3)$ \\ \hline
$E_{7}$ & $(4,2)$ \\ \hline
\end{tabular}%
.
\end{center}

The paper is arranged as follows. In Section 2 we develop fundamental
properties of the Serre spectral sequence of the torus fibration $\pi
:G\rightarrow G/T$, where $T$ is a maximal torus on $G$. The results are
applied in Section 3 to show Theorems A, B and D. Section 4 is devoted to an
exact sequence associated to the cyclic coverings of Lie groups, by which
the proof of Theorem D is reduced to computing with the Chow rings of Lie
groups.

In this paper the cohomologies and spectral sequences are over the ring $%
\mathbb{Z}$ of integers, unless otherwise stated. Given a subset $S$ of a
ring $A$ the notion $\left\langle S\right\rangle $ stands for the ideal
generated by $S$, while $A/\left\langle S\right\rangle $ denotes the
quotient ring. In addition, the elements in a graded ring or algebra are
assumed to be homogeneous.

\bigskip

\noindent \textbf{Remark 1.6.} The problem of computing the cohomology of
Lie groups was raised by E. Cartan in 1929, which was solved only for the
cohomologies with field coefficients, see Reeder \cite{R} and Ka\v{c} \cite%
{K} for accounts about the history. As for the task of the present work
general results on the integral cohomology of Lie groups, such as Theorems A
and B, are requested.

In his problem 15 Hilbert asked for a rigorous foundation of Schubert
calculus. Van der Waerden \cite{Wa} and A. Weil \cite[p.331]{We} attributed
the problem to the determination of the cohomology rings of the flag
manifolds $G/T$. The proofs of Theorems C and D illustrate how Schubert
calculus could be extended to computing with the integral cohomologies of
Lie groups \cite[Remark 6.3]{DZ3}.$\square $

\section{The integral cohomology of Lie groups}

For a Lie group $G$ with a maximal torus $T$ consider the torus fibration

\begin{enumerate}
\item[(2.1)] $T\rightarrow G\overset{\pi }{\rightarrow }G/T$.
\end{enumerate}

\noindent on the group $G$. \textsl{The} \textsl{Borel} \textsl{transgression%
}\textbf{\ }in $\pi $ is the composition

\begin{quote}
$\tau =$ $(\pi ^{\ast })^{-1}\circ \delta :H^{1}(T)\overset{\delta }{%
\rightarrow }H^{2}(G,T)\underset{\cong }{\overset{(\pi ^{\ast })^{-1}}{%
\rightarrow }}H^{2}(G/T)$ (\cite{D1}),
\end{quote}

\noindent where $\delta $ is the connecting homomorphism in the
cohomological exact sequence of the pair $(G,T)$, and where the map $\pi
^{\ast }$ from $H^{2}(G/T)$ to $H^{2}(G,T)$ is always an isomorphism. By the
Leray-Serre theorem \cite[p.135]{Mc} we have that

\bigskip

\noindent \textbf{Lemma 2.1.} \textsl{The second} \textsl{page} \textsl{of
the Serre spectral sequence} $\left\{ E_{r}^{\ast ,\ast }(G),d_{r}\right\} $ 
\textsl{of} \textsl{the fibration} $\pi $ \textsl{is the Koszul complex}

\begin{enumerate}
\item[(2.2)] $E_{2}^{\ast ,\ast }(G)=H^{\ast }(G/T)\otimes H^{\ast }(T)$ 
\textsl{(see} \textsl{\cite[p.259]{Mc})}
\end{enumerate}

\noindent \textsl{on which the differential} $d_{2}$ \textsl{is determined
by the transgression }$\tau $ \textsl{as}

\begin{quote}
\textsl{i)} $d_{2}(x\otimes 1)=0$\textsl{,} $d_{2}(1\otimes t)=\tau
(t)\otimes 1$\textsl{,} $t\in H^{1}(T)$\textsl{, }

\textsl{ii) }$d_{2}(z\cdot z^{\prime })=d_{2}(z)\cdot z^{\prime }+(-1)^{\deg
z}z\cdot d_{2}(z^{\prime })$\textsl{, }$z,z^{\prime }\in E_{2}^{\ast ,\ast
}(G)$.
\end{quote}

\noindent \textsl{In particular, if }$\dim G/T=m$ \textsl{and} $\dim T=n$%
\textsl{, then}

\begin{enumerate}
\item[(2.3)] $E_{2}^{m,n}(G)=H^{m}(G)\otimes
H^{n}(T)=E_{r}^{m,n}(G)=H^{n+m}(G)=\mathbb{Z}$, $r\geq 2$;

\item[(2.4)] $E_{3}^{\ast ,0}(G)=H^{\ast }(G/T)/\left\langle \func{Im}\tau
\right\rangle $\textsl{;}
\end{enumerate}

\noindent \textbf{Proof.} The base manifold $G/T$ of $\pi $ is the complete
flag manifold associated to the Lie group $G$, hence is simply connected.
Formula (2.2), together with properties i) and ii), are standard \cite[p.259]%
{Mc}.

By formula (2.2) we have that

\begin{quote}
$E_{2}^{m,n}(G)=H^{m}(G/T)\otimes H^{n}(T)=\mathbb{Z}$
\end{quote}

\noindent and that $E_{2}^{s,t}(G)=0$ if either $s>m$ or $t>n$. It implies
that any differential $d_{r}$ that acts or lands on the group $%
E_{r}^{m,n}(G) $ must be trivial, confirming the isomorphisms in (2.3).

Finally, by formula i) the differential

\begin{quote}
$d_{2}:E_{2}^{\ast ,1}(G)=H^{\ast }(G/T)\otimes H^{1}(T)\rightarrow
E_{2}^{\ast ,0}(G)=H^{\ast }(G/T)$
\end{quote}

\noindent is $d_{2}(x\otimes t)=x\cup \tau (t)$, $t\in H^{1}(T)$. We get
formula (2.4) from

\begin{quote}
$E_{3}^{\ast ,0}(G)=E_{2}^{\ast ,0}(G)/\func{Im}d_{2}=H^{\ast
}(G/T)/\left\langle \func{Im}\tau \right\rangle $.$\square $
\end{quote}

The Koszul complex $\left\{ E_{2}^{\ast ,\ast }(G),d_{2}\right\} $ has now
been well understood by the following works:

\begin{quote}
a) the base manifold $G/T$ is a flag variety which has a canonical
decomposition into the Schubert cells on $G/T$ \cite{BGG,Ch1};

b) presentation of the ring $H^{\ast }(G/T)$ by special Schubert classes on $%
G/T$ has been completed by Duan and Zhao in \cite{DZ1};

c) with respect to the Schubert basis on $H^{2}(G/T)$ formula of the
transgression $\tau $ is available in \cite[Theorem 2.5]{D1}.
\end{quote}

\noindent Combining these results explicit construction of the ring $%
E_{3}^{\ast ,\ast }(G)$ has been carried out in \cite{D,DZ2}, among which we
shall only need the following two results.

In \cite{G} Grothendieck defined \textsl{the Chow ring of a Lie group} $G$
to be the subring $\mathcal{A}^{\ast }(G):=\pi ^{\ast }(H^{\ast }(G/T))$ of $%
H^{\ast }(G)$, where he has also shown

\begin{quote}
$\mathcal{A}^{\ast }(G)=H^{\ast }(G/T)/\left\langle \func{Im}\tau
\right\rangle $ (\cite[p.21, Rem.2]{G}).
\end{quote}

\noindent Comparing this with formula (2.4) we obtain that

\bigskip

\noindent \textbf{Lemma 2.2.} \textsl{The map }$\pi $\textsl{\ in (2.1)
induces an isomorphism}

\begin{enumerate}
\item[(2.5)] $\mathcal{A}^{\ast }(G)=E_{3}^{\ast ,0}(G)$.
\end{enumerate}

\noindent \textsl{In particular,} \textsl{for the Lie groups }$G$\textsl{\
given in Table 1, the rings }$\mathcal{A}^{\ast }(G)$\textsl{\ and }$%
\mathcal{A}^{\ast }(PG)$\textsl{\ admit the following presentations (in
terms of generators--relations):}

\begin{enumerate}
\item[(2.6)] $\mathcal{A}^{\ast }(SU(n))=\mathbb{Z}$\textsl{,}{\small \ }$%
\mathcal{A}^{\ast }(PSU(n)){\small =}\frac{\mathbb{Z}[\omega ]}{\left\langle
b_{n,r}\omega ^{r},\text{ }1\leq r\leq n\right\rangle }${\small , }$%
b_{n,r}=g.c.d.\{\binom{n}{1},\cdots ,\binom{n}{r}\}$\textsl{;}

\item[(2.7)] $\mathcal{A}^{\ast }(Sp(n))=\mathbb{Z}\QTR{sl}{,}$\textsl{\ }$%
\mathcal{A}^{\ast }(PSp(n))=\frac{\mathbb{Z}[\omega ]}{\left\langle 2\omega
,\omega ^{2^{r+1}}\right\rangle }$\textsl{, where }$n=2^{r}(2k+1)$\textsl{;}

\item[(2.8)] $\mathcal{A}^{\ast }(Spin(2n+1))=\frac{\mathbb{Z[}%
x_{3},x_{5},\cdots ,x_{2\left[ \frac{n+1}{2}\right] -1}]}{\left\langle
2x_{2i-1},\text{ }x_{2i-1}^{k_{i}}\text{; }2\leq i\leq \left[ \frac{n+1}{2}%
\right] \right\rangle }$\textsl{,}

$\mathcal{A}^{\ast }(PSpin(2n+1))=\frac{\mathbb{Z[}x_{1},x_{3},x_{5},\cdots
,x_{2\left[ \frac{n+1}{2}\right] -1}]}{\left\langle 2x_{2i-1},\text{ }%
x_{2i-1}^{k_{i}};\text{ }1\leq i\leq \left[ \frac{n+1}{2}\right]
\right\rangle }$\textsl{;}

\item[(2.9)] $\mathcal{A}^{\ast }(Spin(2n))=\frac{\mathbb{Z[}%
x_{3},x_{5},\cdots ,x_{2\left[ \frac{n}{2}\right] -1}]}{\left\langle
2x_{2i-1},\text{ }x_{2i-1}^{k_{i}}\text{; }2\leq i\leq \left[ \frac{n}{2}%
\right] \right\rangle }$\textsl{\ ;}

$\mathcal{A}^{\ast }(PSpin(2n))=\left\{ 
\begin{tabular}{l}
$\frac{\mathbb{Z[}x_{1},x_{3},x_{5},\cdots ,x_{2\left[ \frac{n}{2}\right]
-1}]}{\left\langle 4x_{1},2x_{1}^{2},x_{1}^{k_{1}},2x_{2i-1},\text{ }%
x_{2i-1}^{k_{i}}\text{; }2\leq i\leq \left[ \frac{n}{2}\right] \right\rangle 
}$ \textsl{if} $n\equiv 1\func{mod}2$\textsl{,} \\ 
$\frac{\mathbb{Z[}\omega ,x_{1},x_{3},x_{5},\cdots ,x_{2\left[ \frac{n}{2}%
\right] -1}]}{\left\langle 2\omega ,\omega ^{2^{h}},2x_{2i-1},\text{ }%
x_{2i-1}^{k_{i}}\text{; }1\leq i\leq \left[ \frac{n}{2}\right] \right\rangle 
}$ \textsl{if }$n=2^{h}(2b+1))$\textsl{,} $h\geq 1$\textsl{;}%
\end{tabular}%
\right. $

\item[(2.10)] $\mathcal{A}^{\ast }(E_{6}){\small =}\frac{\mathbb{Z}%
[x_{3},x_{4}]}{\left\langle 2x_{3},3x_{4},x_{3}^{2},x_{4}^{3}\right\rangle }$%
{\small , }

$\mathcal{A}^{\ast }(PE_{6}){\small =}\frac{\mathbb{Z}[\omega ,x_{3},x_{4}]}{%
\left\langle 3\omega ,2x_{3},3x_{4},(x_{3}+\omega ^{3})^{2},\omega
^{9},x_{4}^{3}\right\rangle }${\small ;}

\item[(2.11)] $\mathcal{A}^{\ast }(E_{7}){\small =}\frac{\mathbb{Z}%
[x_{3},x_{4},x_{5},x_{9}]}{\left\langle 2x_{3},3x_{4},2x_{5},2{x_{9},}%
x_{3}^{2},x_{4}^{3},x_{5}^{2},x_{9}^{2}\right\rangle }${\small , }

$\mathcal{A}^{\ast }(PE_{7}){\small =}\frac{\mathbb{Z}[\omega
,x_{3},x_{4},x_{5},x_{9}]}{\left\langle 2\omega ,\omega
^{2},2x_{3},3x_{4},2x_{5},2{x_{9},}x_{3}^{2},x_{4}^{3},x_{5}^{2},x_{9}^{2}%
\right\rangle }$,
\end{enumerate}

\noindent \textsl{where the generators }$\omega ,$\textsl{\ }$x_{i}$\textsl{%
\ (}$\deg \omega =2$\textsl{, }$\deg x_{k}=2k$\textsl{) are the }$\pi ^{\ast
}$\textsl{--images of certain Schubert classes on }$G/T$\textsl{\ specified
in \cite{DZ1}, and where the power }$k_{i}$\textsl{'s appearing in the
denominator of the quotient (2.7) or (2.8) are respectively}

\begin{quote}
$k_{i}=2^{\left[ \ln \frac{n}{i}\right] +1}$ \textsl{or} $k_{i}=2^{\left[
\ln \frac{n-1}{i}\right] +1}$.
\end{quote}

\noindent \textbf{Proof.} Combining the results mentioned in b) and c) above
it is straightforward to evaluate the Chow rings $\mathcal{A}(G)$ and $%
\mathcal{A}^{\ast }(PG)$ using formula (2.4). As examples we show the
results for $G=PSU(n)$ and $PSpin(2n)$, and refer the remaining cases to 
\cite[Corollary 6.2; formula (6.13)]{DZ1} and \cite[p.i]{M}.

For $G=SU(n)$ the ring $H^{\ast }(G/T)$ has the presentation

\begin{quote}
\noindent $H^{\ast }(G/T)=\frac{\mathbb{Z}[t_{1},t_{2},\cdots ,t_{n}]}{%
\left\langle e_{1},\cdots ,e_{n}\right\rangle }$ (by Borel),
\end{quote}

\noindent where $\deg t_{i}=2$, $e_{r}$ is the $r^{th}$ elementary symmetric
polynomial in the $t_{1},\cdots ,t_{n}$. On the other hand for the group $%
PSU(n)$ we have

\begin{quote}
$\func{Im}\tau =\{t_{1}-t_{2},t_{2}-t_{3},\cdots ,t_{n-1}-t_{n}\}$ by \cite[%
Corollary 3.2]{D1}.
\end{quote}

\noindent From (2.4) we get

\begin{quote}
$\mathcal{A}^{\ast }(PSU(n))=\frac{\mathbb{Z}[t_{1},t_{2},\cdots ,t_{n}]}{%
\left\langle e_{1},\cdots ,e_{n}\right\rangle }\mid _{t_{i}=t_{j}}=\frac{%
\mathbb{Z}[\omega ]}{\left\langle \binom{n}{r}\omega ^{r}\text{, }1\leq
r\leq n\right\rangle }$,
\end{quote}

\noindent where $\omega :=t_{1}=\cdots =t_{n}$. It implies that the order of
the power $\omega ^{r}$ in the ring $\mathcal{A}^{\ast }(PSU(n))$ is
precisely $b_{n,r}=g.c.d.\{\binom{n}{1},\cdots ,\binom{n}{r}\}$, showing
(2.6).

For $G=Spin(2n)$ we have by Marlin \cite[p.20, Theorem 3]{M} that

\begin{quote}
$H^{\ast }(Spin(2n)/T)=\frac{\mathbb{Z}[t_{1},t_{2},\cdots
,t_{n},x_{1},x_{2},\cdots ,x_{n-1}]}{\left\langle \delta _{1},\cdots ,\delta
_{n},\text{ }q_{1},\cdots ,q_{n-1}\right\rangle }$, $\deg t_{i}=2$, $\deg
x_{r}=2r$,
\end{quote}

\noindent where, with $x_{k}=0$ for $k\geq n$ being understood,

\begin{quote}
$\delta _{r}=2x_{r}-e_{r}(t_{1},\cdots ,t_{n})$, $1\leq i\leq n-1$; $\delta
_{n}=e_{n}(t_{1},\cdots ,t_{n})$,

$q_{r}=x_{2r}+2\underset{1\leq i\leq r-1}{\Sigma }%
(-1)x_{i}x_{2r-i}+(-1)^{r}x_{r}^{2}$, $1\leq r\leq n-1$,
\end{quote}

\noindent and where $e_{r}(t_{1},\cdots ,t_{n})$ is the $r^{th}$ elementary
symmetric polynomials in $t_{1},\cdots ,t_{n}$. On the other hand, for the
group $PSpin(2n)$ we have

\begin{quote}
$\func{Im}\tau =\{t_{1}-t_{2},\cdots ,t_{n-1}-t_{n},2t_{n}\}$ by \cite[%
Corollary 3.2]{D1}.
\end{quote}

\noindent The second formula in (2.9) is verified by

\begin{quote}
$\mathcal{A}^{\ast }(PSpin(2n))=H^{\ast }(Spin(2n)/T)/\left\langle \func{Im}%
\tau \right\rangle $

$=H^{\ast }(Spin(2n)/T)\mid _{t_{1}=\cdots =t_{n},2t_{n}=0}$,
\end{quote}

\noindent where $x_{1}=\frac{1}{2}e_{1}$ (by the relation $\delta _{1}$),
and where $\omega :=t_{1}=\cdots =t_{n}$.$\square $

\bigskip

With the product inherited from that on $E_{2}^{\ast ,\ast }(G)$ the third
page $E_{3}^{\ast ,\ast }(G)$ is a bi-graded ring \cite[P.668]{Wh}. In
particular, $E_{3}^{\ast ,1}(G)$ can be considered as a module over the
subring $\mathcal{A}^{\ast }(G)=E_{3}^{\ast ,0}(G)\subset E_{3}^{\ast ,\ast
}(G)$.

\bigskip

\noindent \textbf{Lemma 2.3. }\textsl{For each compact Lie group }$G$\textsl{%
\ with }$I_{G}=\{r_{1},\cdots ,r_{n}\}$\textsl{\ there exist }$n$\textsl{\
elements }$\rho _{k}\in E_{3}^{2r_{k}-2,1}(G)$\textsl{, }$1\leq k\leq n$%
\textsl{, so that}

\textsl{i) as a }$\mathcal{A}^{\ast }(G)$\textsl{--module }$E_{3}^{\ast
,1}(G)$\textsl{\ is spanned additively by }$\left\{ \rho _{1},\cdots ,\rho
_{n}\right\} $\textsl{;}

\textsl{ii) the product }$\rho _{1}\cdots \rho _{n}$\textsl{\ generates the
group }$E_{3}^{m,n}(G)=\mathbb{Z}$\textsl{\ (see (2.3)).}$\square $

\bigskip

\noindent \textbf{Example 2.4.} A classical result of Leray \cite[Example 1.4%
]{D1} states that for the cohomology with real coefficients one has
isomorphism of algebras

\begin{quote}
$H^{\ast }(G;\mathbb{R})=E_{3}^{\ast ,\ast }(G;\mathbb{R})=\Lambda _{\mathbb{%
R}}(y_{1},\cdots ,y_{n})$
\end{quote}

\noindent where the generators $y_{1},\cdots ,y_{n}$ is a basis of the
subspace $E_{3}^{\ast ,1}(G;\mathbb{R})$ constructed from the basic Weyl
invariants of the group $G$. The elements $\rho _{1},\cdots ,\rho _{n}\in
E_{3}^{\ast ,1}(G)$ asserted by Lemma 2.3 may be seen as the integral lifts
of the classes $y_{1},\cdots ,y_{n}$, and will be called \textsl{a basis} of
the $\mathcal{A}^{\ast }(G)$--module $E_{3}^{\ast ,1}(G)$ (for the latter
convenience). Note that $E_{3}^{\ast ,1}(G)$ may fail to be a free module
over the ring $\mathcal{A}^{\ast }(G)$, e.g. see \cite[Example 1.4]{D1} for
the case $G=PSU(8)$.

Without the loss of generalities we can assume in Lemma 2.3 that the Lie
group $G$ is simple. In this case a basis $\left\{ \rho _{1},\cdots ,\rho
_{n}\right\} $ of the module $E_{3}^{\ast ,1}(G)$ has been constructed
uniformly for all simply connected $G$ in \cite{DZ2}, and for the non-simply
connected $G$ in \cite{D,D1}. In this paper knowing the degrees of these
elements is sufficient for our purpose.

It is worth to mention that, in the course to compute the integral
cohomology of the spinor group $Spin(n)$ Pittie \cite{P} has constructed a
basis $\left\{ \rho _{1},\rho _{2},\cdots \right\} $ of the module $%
E_{3}^{\ast ,1}(G)$ for the classical groups $G=U(n),SO(n)$ and $Spin(n)$.$%
\square $

\bigskip

To see the implication of Lemma 2.3 let $F^{p}$ be the filtration on the
integral cohomology $H^{\ast }(G)$ defined by the map $\pi $. That is (\cite[%
P.146]{Mc})

\begin{center}
$0=F^{r+1}(H^{r}(G))\subseteq F^{r}(H^{r}(G))\subseteq \cdots \subseteq
F^{0}(H^{r}(G))=H^{r}(G)$
\end{center}

\noindent with

\begin{quote}
$E_{\infty }^{p,q}(G)=F^{p}(H^{p+q}(G))/F^{p+1}(H^{p+q}(G))$.
\end{quote}

\noindent The fact $H^{2k+1}(G/T)=0$ due to Bott and Samelson \cite{BS}
implies that

\begin{quote}
a) $E_{r}^{2k+1,\ast }=0$ for $k\geq 0$; \quad b) $%
E_{3}^{4s,2}=E_{4}^{4s,2}=\cdots =E_{\infty }^{4s,2}$.
\end{quote}

\noindent From $F^{2s+1}(H^{2s+1}(G))=F^{2s+2}(H^{2s+1}(G))=0$ by a) one
finds that

\begin{quote}
$E_{\infty }^{2s,1}(G)=F^{2s}(H^{2s+1}(G))\subset H^{2s+1}(G)$.
\end{quote}

\noindent Combining this with the routine relation $d_{r}(E_{r}^{\ast ,1})=0$
for $r\geq 3$ yields the composition

\begin{enumerate}
\item[(2.12)] $\kappa :E_{3}^{\ast ,1}(G)\rightarrow E_{4}^{\ast
,1}(G)\rightarrow \cdots \rightarrow E_{\infty }^{\ast ,1}(G)\subset H^{\ast
}(G)$
\end{enumerate}

\noindent that interprets elements of $E_{3}^{\ast ,1}$ directly as
cohomology classes of $G$. Carrying on the results of Lemma 2.3 we obtain
the following characterization of the integral cohomology $H^{\ast }(G)$ in
term of its free part and torsion ideal.

\bigskip

\noindent \textbf{Lemma 2.5. }\textsl{The integral cohomology of }$G$ 
\textsl{has the additive presentation}

\begin{enumerate}
\item[(2.13)] $H^{\ast }(G)=\Delta (\xi _{1},\cdots ,\xi _{n})\oplus \tau
(G) $\textsl{\ with }$\xi _{i}:=\kappa (\rho _{i})\in H^{\ast }(G)$\textsl{,}
\end{enumerate}

\noindent \textsl{where }$\Delta (\xi _{1},\cdots ,\xi _{n})$ \textsl{%
denotes the free} $\mathbb{Z}$\textsl{--module with the basis}

\begin{quote}
$\Phi :=\{1,\xi _{I}=\underset{i\in I}{\cup }\xi _{i}\in H^{\ast }(G)\mid
I\subseteq \{1,\cdots ,n\}\}$.
\end{quote}

\noindent \textbf{Proof.} Let $\left\{ E_{r}^{\ast ,\ast }(G;\mathbb{R}%
),d_{r}\right\} $\ be the Serre spectral sequence of $\pi $ with real
coefficients. According to Leray \cite{L} the algebra $E_{3}^{\ast ,\ast }(G,%
\mathbb{R})$ is generated multiplicatively by its subspace $E_{3}^{\ast
,1}(G,\mathbb{R})$, while the map $\kappa $ in (2.12) induces

\begin{quote}
i) an isomorphism $E_{3}^{\ast ,1}(G,\mathbb{R})\cong P(G;\mathbb{R})$ of
vector spaces, and

ii) an isomorphism $E_{3}^{\ast ,\ast }(G,\mathbb{R})\cong H^{\ast }(G;%
\mathbb{R)}$ of algebras,
\end{quote}

\noindent respectively. Assume that $\left\{ \rho _{1},\cdots ,\rho
_{n}\right\} $ is a basis of the $\mathcal{A}^{\ast }(G)$--module $%
E_{3}^{\ast ,1}(G)$, and that $\dim G/T=m$, $\dim T=n$. In addition to $\xi
_{i}:=\kappa (\rho _{i})$ we put

\begin{quote}
$\xi _{i}^{R}:=\xi _{i}\otimes 1\in $ $H^{\ast }(G;\mathbb{R})=H^{\ast
}(G)\otimes \mathbb{R}$, $1\leq i\leq n$.
\end{quote}

\noindent Since $\mathcal{A}^{\ast }(G)\otimes \mathbb{R}=\mathbb{R}$ the
space $E_{3}^{\ast ,1}(G;\mathbb{R})=E_{3}^{\ast ,1}(G)\otimes \mathbb{R}$%
\textsl{\ }has the basis\textsl{\ }$\left\{ \rho _{1}\otimes 1,\cdots ,\rho
_{n}\otimes 1\right\} $ by i) of Lemma 2.3. By the isomorphisms i) and ii)
above

\begin{quote}
$H^{\ast }(G;\mathbb{R)=}\Lambda _{\mathbb{R}}(\xi _{1}^{R},\cdots ,\xi
_{n}^{R})$ with $\xi _{i}^{R}\in P(G;\mathbb{R})$,
\end{quote}

\noindent implying that the set $\Phi $ of $2^{n}$ monomials is linearly
independent in $H^{\ast }(G)$. It remains to show that the set $\Phi $ spans
a direct summand of $H^{\ast }(G)$.

Assume, on the contrary, that there exist a monomial $\xi _{I}\in \Phi $, an
integral class $\varsigma \in H^{\ast }(G)$, as well as some integer $a>1$,
so that a relation of the form $\xi _{I}=a\cdot \varsigma $ holds in $%
H^{\ast }(G)$. Multiplying both sides by $\xi _{J}$ with $J$ the complement
of $I\subseteq \{1,\cdots ,n\}$ yields that

\begin{quote}
$\xi _{1}\cup \cdots \cup \xi _{n}=(-1)^{r}a\cdot (\varsigma \cup \xi _{J})$
(for some $r\in \{0,1\}$).
\end{quote}

\noindent However, in view of the identification (2.3) the map $\kappa $ in
(2.12) transforms the generator $\rho _{1}\cdots \rho _{n}$ of $%
E_{3}^{m,n}(G)=\mathbb{Z}$ to the generator $\xi _{1}\cup \cdots \cup \xi
_{n}$ of $H^{n+m}(G)=\mathbb{Z}$ (by $\xi _{i}:=\kappa (\rho _{i})$). The
proof is completed by this contradiction.$\square $

\section{Proofs of Theorems A, B and C}

Let $\left\{ \rho _{1},\cdots ,\rho _{n}\right\} $ be a basis of the $%
\mathcal{A}^{\ast }(G)$--module\textsl{\ }$E_{3}^{\ast ,1}(G)$. With $\xi
_{i}=\kappa (\rho _{i})$ and $\xi _{i}^{R}=\xi _{i}\otimes 1$ we have by the
proof of Lemma 2.5 that

\begin{enumerate}
\item[(3.1)] $H^{\ast }(G;\mathbb{R})=\Lambda _{\mathbb{R}}(\xi
_{1}^{R},\cdots ,\xi _{n}^{R})$, where $\xi _{i}^{R}\in P(G;\mathbb{R})$%
\textsl{.}
\end{enumerate}

\noindent \textbf{Proof of Theorem A.} Let $q:H^{\ast }(G)\rightarrow 
\mathcal{F}(G)=H^{\ast }(G)/\tau (G)$ be the quotient map and put $%
x_{i}:=q(\xi _{i})$. Then $\mathcal{F}(G)=\Delta (x_{1},\cdots ,x_{n})$ by
(2.13), where

\begin{quote}
$\deg x_{i}=\deg \xi _{i}=2r_{i}-1,1\leq i\leq n$, $I_{G}=\{r_{1},\cdots
,r_{n}\}$.
\end{quote}

\noindent In addition, with $\deg \xi _{i}\equiv 1\func{mod}2$ we get $%
x_{i}^{2}=0$ from $\xi _{i}^{2}\in \tau (G)$, showing

\begin{enumerate}
\item[(3.2)] $\mathcal{F}(G)=\Lambda (x_{1},\cdots ,x_{n})$ (i.e. the
formula i) of Theorem A).
\end{enumerate}

For the co--product $\mu ^{\ast }$ on $H^{\ast }(G)$ we can assume, in
general, that

\begin{quote}
$\mu ^{\ast }(\xi _{i})=\xi _{i}\otimes 1+1\otimes \xi _{i}+z_{i}+d_{i}$, $%
1\leq i\leq n$,
\end{quote}

\noindent where $z_{i}=\Sigma a_{j}\otimes b_{j}\in \mathcal{F}^{+}(G)\times 
\mathcal{F}^{+}(G)$ is a mixed term, $d_{i}\in \tau (G\times G)$. With $\xi
_{i}^{R}\in $ $P(G;\mathbb{R})$ by (3.1) there must be $z_{i}=0$. Thus,
applying $q\times q$ to this equality verifying the assertion ii) of Theorem
A

\begin{quote}
$\mu ^{\#}(x_{i})=x_{i}\otimes 1+1\otimes x_{i}$.
\end{quote}

Finally, suppose that the group $G$\ is simple, and that in addition to
(3.2) one has $\mathcal{F}(G)=\Lambda (x_{1}^{\prime },\cdots ,x_{n}^{\prime
})$ with

\begin{quote}
$\deg x_{i}=\deg x_{i}^{\prime }$ and $\mu ^{\#}(x_{i}^{\prime
})=x_{i}^{\prime }\otimes 1+1\otimes x_{i}^{\prime }$\textsl{, }$1\leq i\leq
n$.
\end{quote}

\noindent Let $D\subset \mathcal{F}(G)$ be the subring consists of the
decomposable elements in the positive degrees. Since both $\left\{
x_{1}^{\prime },\cdots ,x_{n}^{\prime }\right\} $ and $\left\{ x_{1},\cdots
,x_{n}\right\} $ are basis of the free $\mathbb{Z}$--module $\mathcal{F}%
(G)/D $, and since $2=r_{1}<\cdots <r_{n}$, the transition functions on $%
\mathcal{F}(G)$ between these two sets of generators take the forms

\begin{quote}
$x_{i}^{\prime }=\pm x_{i}+d_{i}$ with $d_{i}\in D$, $1\leq i\leq n$.
\end{quote}

\noindent Applying the co--product $\mu ^{\#}$ to both sides, and using the
relations $\mu ^{\#}(y)=y\otimes 1+1\otimes y$ for $y=x_{i}$ or $%
x_{i}^{\prime }$ to simplify the resulting equation, yields that

\begin{quote}
$d_{i}\otimes 1+1\otimes d_{i}=\mu ^{\#}(d_{i})$.
\end{quote}

\noindent With $d_{i}\in D$ this implies $d_{i}=0$, completing the proof of
Theorem A.$\square $

\bigskip

\noindent \textbf{Proof of Theorem B.} Let $G$ be a simple Lie group with
maximal torus $T$ and non--trivial center $\mathcal{Z}(G)$. Since the
covering $c:G\rightarrow PG=G/\mathcal{Z}(G)$ carries $T$ to the maximal
torus $T^{\prime }:=c(T)$ on $PG$, it can be viewed as a bundle map over the
identity of $G/T$

\begin{quote}
$%
\begin{array}{ccc}
T & \overset{c^{\prime }}{\rightarrow } & T^{\prime } \\ 
\cap \quad &  & \cap \quad \\ 
G & \overset{c}{\rightarrow } & PG \\ 
\pi \downarrow \quad &  & \pi ^{\prime }\downarrow \quad \\ 
G/T & = & PG/T^{\prime }%
\end{array}%
$ (note that $G/T=PG/T^{\prime }$),
\end{quote}

\noindent hence induces a map $c^{\ast }:E_{2}^{\ast ,\ast }(PG)\rightarrow
E_{2}^{\ast ,\ast }(G)$ of Koszul complexes. Assume by Lemma 2.3 that $%
\left\{ \rho _{1}^{\prime },\cdots ,\rho _{n}^{\prime }\right\} $ is a basis
of the $\mathcal{A}(PG)$--module $E_{3}^{\ast ,1}(PG)$, and that $\left\{
\rho _{1},\cdots ,\rho _{n}\right\} $ is a basis of the $\mathcal{A}^{\ast
}(G)$--module $E_{3}^{\ast ,1}(G)$. Then, by the proof of Theorem A

\begin{quote}
$\mathcal{F}(PG)=\Lambda (y_{1},\cdots ,y_{n})$ with $y_{i}:=q\circ \kappa
(\rho _{i}^{\prime })\in \mathcal{F}(PG)$,

$\mathcal{F}(G)=\Lambda (x_{1},\cdots ,x_{n})$ with $x_{i}:=q\circ \kappa
(\rho _{i})\in \mathcal{F}(G)$.
\end{quote}

\noindent On the other hand, for the degree reason $\deg \rho _{i}=\deg \rho
_{i}^{\prime }$ we can assume by i) and Lemma 2.4 that

\begin{enumerate}
\item[(3.3)] $c^{\ast }(\rho _{i}^{\prime })=a_{i}\rho _{i}+b_{i-1}\rho
_{i-1}+\cdots +b_{1}\rho _{1}$, $a_{i}\in \mathbb{Z}$, $b_{j}\in \mathcal{A}%
^{+}(G)$,
\end{enumerate}

\noindent where $\mathcal{A}^{+}(G)$ denotes the subring of $\mathcal{A}%
^{\ast }(G)$ consisting of elements in the positive degrees. Since the ring $%
\mathcal{A}^{+}(G)$ is always finite (3.3) implies that $c^{\#}(y_{i})=a_{i}%
\cdot x_{i}$ on the quotient $\mathcal{F}(G)$, where we can assume $%
a_{i}\geq 0$ by modifying the sign of the generator $x_{i}$ whenever is
necessary, to get the assertion i) of Theorem B.

Finally, applying $c^{\ast }$ to the generator $\rho _{1}^{\prime }\cdots
\rho _{n}^{\prime }$ of $E_{3}^{m,n}(PG)=H^{m+n}(PG)=\mathbb{Z}$ (see Lemma
2.3 and (2.3)), and noting that the partial sum $b_{i-1}\rho _{i-1}+\cdots
+b_{1}\rho _{1}$ in (3.3) is of finite order, we get

\begin{enumerate}
\item[(3.4)] $c^{\ast }(\rho _{1}^{\prime }\cdots \rho _{n}^{\prime
})=(a_{1}\cdots a_{n})\rho _{1}\cdots \rho _{n}$ on $%
E_{3}^{m,n}(G)=H^{m+n}(G)=\mathbb{Z}$.
\end{enumerate}

\noindent Since the mapping degree of $c$ is the order $\left\vert \mathcal{Z%
}(G)\right\vert $ of the center, and since the product $\rho _{1}\cdots \rho
_{n}$ generates the group $E_{3}^{m,n}(G)=H^{m+n}(G)=\mathbb{Z}$ by Lemma
2.3, we get from (3.4) that $a_{1}\cdots a_{n}=\left\vert \mathcal{Z}%
(G)\right\vert $, completing the proof of Theorem B.$\square $

\bigskip

\noindent \textbf{Proof of Theorem C. }Since the second homotopy group of $%
PG $ is trivial, the homotopy exact sequence of $\pi ^{\prime }$ contains
the free resolution of the group $\pi _{1}(PG)$

\begin{quote}
$0\rightarrow \pi _{2}(G/T)\rightarrow \pi _{1}(T^{\prime })\overset{j_{\ast
}}{\rightarrow }\pi _{1}(PG)(=\mathcal{Z}(G))\rightarrow 0$,
\end{quote}

\noindent where $j$ is the inclusion of the maximal torus. Applying the
co--functor $Hom(,\mathbb{Z})$ to this sequence, and using the Hurewicz
isomorphisms

\begin{quote}
$\pi _{2}(G/T)=H_{2}(G/T)$, $\pi _{1}(T^{\prime })=H_{1}(T^{\prime })$, $\pi
_{1}(PG)=H_{1}(PG)$
\end{quote}

\noindent to substitute the relevant groups, one obtains the exact sequence

\begin{enumerate}
\item[(3.6)] $0\rightarrow H^{1}(PG)\overset{j^{\ast }}{\rightarrow }%
H^{1}(T^{\prime })\overset{\tau ^{\prime }}{\rightarrow }H^{2}(G/T)\overset{%
\pi ^{\prime \ast }}{\rightarrow }TorH^{2}(PG)=\mathcal{Z}(G)\rightarrow 0$
\end{enumerate}

\noindent in cohomologies, where $\tau ^{\prime }$ is the transgression in $%
\pi ^{\prime }$. It follows that

\begin{quote}
$H^{2}(G/T)/\func{Im}\tau ^{\prime }\cong H^{2}(PG)=TorH^{2}(PG)=\mathcal{Z}%
(G)$,
\end{quote}

\noindent where the first equality follows from $H^{2}(PG)\otimes \mathbb{R}%
=0$ by Theorem 1.1. With $\mathcal{A}^{2}(PG)=H^{2}(G/T)/\func{Im}\tau
^{\prime }$ this implies that the inclusion $\mathcal{A}(PG)\subseteq
H^{\ast }(PG)$ restricts to an isomorphism in degree $2$. Therefore, the
subring $J(PG)$ of $H^{\ast }(PG)$ generated by $H^{2}(PG)$ agrees with the
subring of $\mathcal{A}^{\ast }(PG)$ generated by $\mathcal{A}^{2}(PG)$. In
particular, one reads out $J(PG)$ from the formulae of the ring $\mathcal{A}%
^{\ast }(PG)$ presented in Lemma 2.2, showing Theorem C.$\square $

\section{An exact sequence for cyclic coverings}

As a covering $c:G\rightarrow G^{\prime }$ on Lie groups is a group
homomorphism we have $\ker c\subseteq \mathcal{Z}(G)$. The covering $c$ is
called \textsl{cyclic} if $\ker c$ is a cyclic subgroup. We note that

a) if $G$ is one of the simply connected Lie groups in Table 1 with $G\neq
Spin(2n)$, then the covering $G\rightarrow PG$ is always cyclic;

b) if $G=Spin(2n)$ the covering $G\rightarrow PG$ can be decomposed into the
composition $Spin(2n)\overset{c_{1}}{\rightarrow }SO(2n)\overset{c_{2}}{%
\rightarrow }PSpin(2n)$ of two cyclic ones, both with order $2$.

\noindent In addition, for a cyclic covering $G\rightarrow G^{\prime }$
between simple Lie groups, the multi--degree $\mathcal{D}(G,G^{\prime })$ is
also defined by the proof of Theorem B. Summarizing, to show Theorem D it
suffices to compute the invariant $\mathcal{D}(G,G^{\prime })$ for the
cyclic coverings. In this section we establish an\textsl{\ }exact sequence,%
\textsl{\ }that reduces this task to the calculation with the Chow ring $%
\mathcal{A}^{\ast }(G^{\prime })$ of the group $G^{\prime }$.

Assume that $c:G\rightarrow G^{\prime }$ is a cyclic covering on simple Lie
groups. \textsl{The central extension of} $c$ is the principle $U(1)$%
--fibration over the group $G^{\prime }$

\begin{enumerate}
\item[(4.1)] $0\rightarrow U(1)\rightarrow \widetilde{G}:=G\times _{\ker
c}U(1)\overset{C}{\rightarrow }G^{\prime }\rightarrow 0$,
\end{enumerate}

\noindent where $\ker c$ acts on the cycle $U(1)$ as the anti--clockwise
rotation through the angle $2\pi /\left\vert \ker c\right\vert $, and where $%
\widetilde{G}$ is furnished with the obvious group structure. Moreover, fix
once for all a maximal torus $T$ on $G$, $\dim T=n$. Then both

\begin{quote}
$\widetilde{T}:=T\times _{\ker c}U(1)\subset \widetilde{G}$ and $T^{\prime
}:=c(T)\subset G^{\prime }$
\end{quote}

\noindent are respectively maximal torus of the corresponding Lie groups,
while the $U(1)$--fibration $C$ can be regarded as a bundle map between two
torus fibrations:

\begin{enumerate}
\item[(4.2)] $%
\begin{array}{ccccc}
U(1) & \hookrightarrow & \widetilde{T} & \overset{C^{\prime }}{\rightarrow }
& T^{\prime } \\ 
\parallel &  & \cap \quad &  & \cap \quad \\ 
U(1) & \hookrightarrow & \widetilde{G} & \overset{C}{\rightarrow } & 
G^{\prime } \\ 
&  & \widetilde{\pi }\downarrow \quad &  & \pi ^{\prime }\downarrow \quad \\ 
&  & G/T & = & G/T%
\end{array}%
$ (see (2.1)).
\end{enumerate}

\noindent Since both $\widetilde{T}$ and $T^{\prime }$ are torus groups, the
restriction $C^{\prime }$ of $C$ on $\widetilde{T}$ is splittable, implying
that the group $H^{1}(\widetilde{T})$ admits a basis $\{t_{0},t_{1},\cdots
,t_{n}\}$ so that the induced map $C^{\prime \ast }$ carries the ring $%
H^{\ast }(T^{\prime })$ isomorphically to the subring $\Lambda (t_{1},\cdots
,t_{n})$ of $H^{\ast }(\widetilde{T})=\Lambda (t_{0},t_{1},\cdots ,t_{n})$.
It follows that

\bigskip

\noindent \textbf{Lemma 4.1.} \textsl{For a cyclic covering }$c:G\rightarrow
G^{\prime }$\textsl{\ of simple Lie groups} \textsl{we have}

\begin{quote}
$E_{2}^{\ast ,\ast }(\widetilde{G})=H^{\ast }(G/T)\otimes \Lambda
(t_{0},t_{1},\cdots ,t_{n})$\textsl{; }

$E_{2}^{\ast ,\ast }(G^{\prime })=H^{\ast }(G/T)\otimes \Lambda
(t_{1},\cdots ,t_{n})$,
\end{quote}

\noindent \textsl{on which }$C$ \textsl{induces a map} $C^{\ast }:$\textsl{\ 
}$E_{2}^{\ast ,\ast }(G^{\prime })\rightarrow E_{2}^{\ast ,\ast }(\widetilde{%
G})$ \textsl{of Koszul complexes with}

\begin{quote}
\textsl{i)} $C^{\ast }(z\otimes t)=z\otimes t$\textsl{, where }$z\in H^{\ast
}(G/T)$\textsl{, }$t\in H^{\ast }(T^{\prime });$

\textsl{ii)} \textsl{the transgression }$\tau ^{\prime }$ \textsl{of }$\pi
^{\prime }$\textsl{\ is the restriction of the transgression }$\widetilde{%
\tau }$ \textsl{of} $\widetilde{\pi }$ \textsl{to the subgroup} $%
H^{1}(T^{\prime })\subset H^{1}(\widetilde{T})$.

\textsl{iii)} \textsl{the} \textsl{Euler class of the }$U(1)$\textsl{%
--fibration }$C$ \textsl{over} $G^{\prime }$ \textsl{is }

$\qquad \omega =\left[ \widetilde{\tau }(t_{0})\otimes 1\right] \in $\textsl{%
\ }$E_{3}^{2,0}(G^{\prime })=H^{2}(G^{\prime })$\textsl{,}
\end{quote}

\noindent \textsl{where }$[x]$\textsl{\ denotes the cohomology class of a }$%
d_{2}^{\prime }$\textsl{--cocycle }$x\in E_{2}^{2,k}(G^{\prime })$\textsl{,
and where the identification }$E_{3}^{2,0}(G^{\prime })=H^{2}(G^{\prime })$ 
\textsl{follows from the proof of Theorem C.}$\square $

\bigskip

By Lemma 4.1 the map $C^{\ast }$\textsl{\ }fits into the short exact sequence

\begin{enumerate}
\item[(4.3)] $0\rightarrow E_{2}^{\ast ,k}(G^{\prime })\overset{C^{\ast }}{%
\rightarrow }E_{2}^{\ast ,k}(\widetilde{G})\overset{\theta }{\rightarrow }%
E_{2}^{\ast ,k-1}(G^{\prime })\rightarrow 0$,
\end{enumerate}

\noindent of Koszul complexes, where the map $\theta $ is evaluated by the
simple rule:

\begin{enumerate}
\item[(4.4)] $\theta (x\otimes y)=x\otimes y_{2}$ if $y=y_{1}+t_{0}\cdot
y_{2}$ with $y_{1},y_{2}\in H^{\ast }(T^{\prime })$.
\end{enumerate}

\noindent Moreover, with $rankH^{2}(G/T)=n$ and $\dim \widetilde{T}=n+1$,
the transgression $\widetilde{\tau }$ satisfies that $\ker \widetilde{\tau }=%
\mathbb{Z\subseteq }H^{1}(\widetilde{T})$. Taking a generator $s\in $ $\ker 
\widetilde{\tau }$ and noting that $1\otimes s\in E_{2}^{0,1}(\widetilde{G})$
is $\widetilde{d}_{2}$\textsl{--}closed by i) of Lemma 2.1, we get the
cohomology class

\begin{quote}
$\rho _{0}=[1\otimes s]\in E_{3}^{0,1}(\widetilde{G})$.
\end{quote}

\noindent \textbf{Theorem 4.2.} \textsl{For a cyclic covering }$%
c:G\rightarrow G^{\prime }$ \textsl{one has the exact sequence}

\begin{enumerate}
\item[(4.5)] $0\rightarrow E_{3}^{\ast ,\ast }(G^{\prime })/\left\langle
\omega \right\rangle \overset{{\small C}^{\ast }}{{\small \rightarrow }}%
E_{3}^{\ast ,\ast }(\widetilde{G})\overset{\theta }{\rightarrow }E_{3}^{\ast
,\ast }(G^{\prime })\overset{{\small \omega }}{{\small \rightarrow }}\omega
\cdot E_{3}^{\ast ,\ast }(G^{\prime })\rightarrow 0$,
\end{enumerate}

\noindent \textsl{in which}

\begin{enumerate}
\item[(4.6)] $E_{3}^{\ast ,\ast }(\widetilde{G})=E_{3}^{\ast ,\ast
}(G)\otimes \Lambda (\rho _{0})$\textsl{,}
\end{enumerate}

\noindent \textsl{where }$\omega $ \textsl{is the Euler class of the central
extension }$C$\textsl{\ of }$c$\textsl{, and where} \textsl{for any }$z\in
E_{3}^{\ast ,\ast }(G^{\prime })$, $y\in E_{3}^{\ast ,\ast }(\widetilde{G})$%
\textsl{,}

\begin{quote}
\textsl{i)} $\theta (C^{\ast }(z)\cdot y)=z\cdot \theta (y)$\textsl{;}

\textsl{ii)} $\theta (\rho _{0})=\left\vert \ker c\right\vert \in
E_{3}^{0,0}(G^{\prime })\mathbb{=Z}$\textsl{;}

\textsl{iii)} $\omega (z)=\omega \cdot z$\textsl{;}

\textsl{iv)} ${\small C}^{\ast }(z)\equiv c^{\ast }(z)\otimes 1\func{mod}%
\rho _{0}$ \textsl{(with respect to (4.6)).}
\end{quote}

\noindent \textbf{Proof. }The short exact sequence (4.3) of Koszul complexes
yields the long exact sequence in cohomologies

\begin{quote}
${\small \cdots \rightarrow }E_{3}^{\ast ,r}(G^{\prime })\overset{{\small C}%
^{\ast }}{{\small \rightarrow }}E_{3}^{\ast ,r}(\widetilde{G})\overset{%
\overline{\theta }}{{\small \rightarrow }}E_{3}^{\ast ,r-1}(G^{\prime })%
\overset{{\small \omega }}{{\small \rightarrow }}E_{3}^{\ast ,r-1}(G^{\prime
})${\tiny \ }$\overset{{\small C}^{\ast }}{{\small \rightarrow }}{\small %
\cdots }$
\end{quote}

\noindent for which properties i), ii) and iii) follows easily from Lemma
4.1. Since $\ker C^{\ast }=$ $\omega {\small \cdot }E_{3}^{\ast ,\ast
}(G^{\prime })\subset E_{3}^{\ast ,\ast }(G^{\prime })$ (by the exactness)
we obtain (4.5) from

\begin{quote}
$\func{Im}\omega =\omega \cdot E_{3}^{\ast ,\ast }(G^{\prime })$ and $co\ker
\omega =E_{3}^{\ast ,\ast }(G^{\prime })/\func{Im}\omega =E_{3}^{\ast ,\ast
}(G^{\prime })/\left\langle \omega \right\rangle $.
\end{quote}

\noindent It remains to establish the decomposition (4.6), together with the
relation iv).

In addition to (4.2) the group $\widetilde{G}$ also fits into the fibrations

\begin{enumerate}
\item[(4.7)] $%
\begin{array}{ccccccc}
0\rightarrow & T & \hookrightarrow & \widetilde{T} & \overset{g^{\prime }}{%
\rightarrow } & U(1) & \rightarrow 0 \\ 
& \cap \quad &  & \cap \quad &  & \parallel &  \\ 
0\rightarrow & G & \overset{i}{\hookrightarrow } & \widetilde{G} & \overset{g%
}{\rightarrow } & U(1) & \rightarrow 0 \\ 
& \pi \downarrow &  & \widetilde{\pi }\downarrow \quad &  & \quad &  \\ 
& G/T & = & G/T &  &  & 
\end{array}%
$,
\end{enumerate}

\noindent where $g$ is the quotient of the projection $G\times
U(1)\rightarrow U(1)$ by $\ker c$; $g^{\prime }$ is the restriction of $g$
on $\widetilde{T}$; and where the upper two rows are short exact sequences
of Lie groups. It follows from Lemma 2.1 that $E_{2}^{\ast ,\ast }(%
\widetilde{G})=E_{2}^{\ast ,\ast }(G)\otimes \Lambda (s)$. With $\rho
_{0}=[1\otimes s]$ by our convention we obtain (4.6) from the K\"{u}nnth
formula. Finally, the relation iv) is transparent, since the inclusion $i$
in (4.7) identifies $G$ with the normal subgroup $q^{-1}(1)$ of $\widetilde{G%
}$ that satisfies the relation $C\circ i=c$.$\square $

\bigskip

Recall that a $U(1)$--fibration $E\rightarrow X$ over a CW--complex $X$ is
classified by its Euler class $\omega \in H^{2}(X)$. The following result,
fairly transparent in the context of \cite{BH}, provides a geometric
interpretation of the generators of the rings $J(PG)$ given in Theorem C.

\bigskip

\noindent \textbf{Lemma 4.3. }\textsl{Let }$G$\textsl{\ be a simply
connected simple Lie group given in Table 1.}

\textsl{i) if }$G\neq Spin(2^{t}(2b+1))$\textsl{\ with }$t\geq 2$\textsl{,
the Euler class} \textsl{of the central extension\ of the cyclic covering }$%
c:G\rightarrow PG$\textsl{\ is the generator }$\omega \in J(PG)$\textsl{;}

\textsl{ii) if }$G=Spin(2n)$\textsl{\ with }$n\equiv 0\func{mod}2$\textsl{,} 
\textsl{the Euler class} \textsl{of the central extension of the 2 sheets
covering }$c:SO(2n)\rightarrow PSpin(2n)$\textsl{\ is the generator }$\omega
_{1}\in J(PSpin(2n))$\textsl{.}$\square $

\bigskip

It is crucial to note in (4.5) that, with respect to the bi--gradation on $%
E_{3}^{s,t}$ imposed by the "\textsl{base degrees }$s$" and "\textsl{fiber
degrees }$t$", the map $C^{\ast }$ preserves both the base and fiber
degrees; the map $\theta $ preserves the base degrees, but reduces the fiber
degrees by $1$; and that the map $\omega $ increases the base degrees by $2$%
, and preserves the fiber degrees. In particular, for each $k\geq 0$ one has
by (4.5) the exact sequences with four terms

\begin{enumerate}
\item[(4.8)] $0\rightarrow E_{3}^{2k,1}(PG)/\left\langle \omega
\right\rangle \overset{C^{\ast }}{\rightarrow }E_{3}^{2k,1}(\widetilde{G})%
\overset{\theta }{\rightarrow }\mathcal{A}^{2k}(PG)\overset{\omega }{%
\rightarrow }\omega \cdot \mathcal{A}^{2k}(PG)\rightarrow 0$
\end{enumerate}

\noindent where the groups $E_{3}^{2k,0}(PG)=\mathcal{A}^{2k}(PG)$ and $%
E_{3}^{2k,0}(\widetilde{G})=\mathcal{A}^{2k}(G)$ have been decided by Lemma
2.2. Furthermore, if $\left\{ \rho _{1}^{\prime },\cdots ,\rho _{n}^{\prime
}\right\} $ is a basis of the $\mathcal{A}^{\ast }(PG)$--module $E_{3}^{\ast
,1}(PG)$, and $\left\{ \rho _{1},\cdots ,\rho _{n}\right\} $ is a basis of
the $\mathcal{A}^{\ast }(G)$--module $E_{3}^{\ast ,1}(G)$, then, with $%
\mathcal{A}(PG)/\left\langle \omega \right\rangle =\mathcal{A}(G)$ by Lemma
2.2,

\begin{enumerate}
\item[(4.9)] $E_{3}^{\ast ,1}(\widetilde{G})$ is a $\mathcal{A}^{\ast }(G)$%
--module with basis $\left\{ \rho _{0},\rho _{1},\cdots ,\rho _{n}\right\} $
by (4.6);

$E_{3}^{\ast ,1}(PG)/\left\langle \omega \right\rangle $ is a $\mathcal{A}%
^{\ast }(G)$--module with basis $\left\{ \rho _{1}^{\prime },\cdots ,\rho
_{n}^{\prime }\right\} $.
\end{enumerate}

\noindent \textbf{Proof of Theorem D. }Let $G$ be a simply connected Lie
group given in Table 1, and assume that $D(G,PG)=\{a_{1},\cdots ,a_{n}\}$.
By ii) of Theorem B we have

\begin{enumerate}
\item[(4.10)] $a_{1}\cdot \cdots \cdot a_{n}=\left\vert \mathcal{Z}%
(G)\right\vert $.
\end{enumerate}

\noindent In addition, by the naturality of the map $\kappa $ in (2.12) with
respect to bundle maps, we have the commutative diagram

\begin{enumerate}
\item[(4.11)] 
\begin{tabular}{lllll}
$E_{3}^{2k,1}(PG)$ & $\overset{\kappa }{\rightarrow }$ & $H^{\ast }(PG)$ & $%
\overset{q}{\rightarrow }$ & $\mathcal{F}(PG)$ \\ 
$C^{\ast }\downarrow $ &  & $C^{\ast }\downarrow $ &  & $C^{\#}\downarrow $
\\ 
$E_{3}^{2k,1}(\widetilde{G})$ & $\overset{\kappa }{\rightarrow }$ & $H^{\ast
}(\widetilde{G})$ & $\overset{q}{\rightarrow }$ & $\mathcal{F}(\widetilde{G}%
) $%
\end{tabular}
\end{enumerate}

\noindent allowing us to calculate the \textsl{eigenvalues} $a_{1},\cdots
,a_{n}$ of $c^{\#}$ (see i) of Theorem B) by computing with the action of $%
C^{\ast }$ on $E_{3}^{2k,1}(PG)$. Granted with the relations (4.8), (4.9)
and (4.10), together with the diagram (4.11), the proof of Theorem D will be
given in the following order

\begin{quote}
$G=SU(n)$, $Sp(n)$, $E_{7}$, $E_{6}$, $Spin(2n+1)$ and $Spin(2n)$,
\end{quote}

\noindent where $\mathbb{Z}_{m}\{x\}$ (resp. $\mathbb{Z}\{x\}$) denotes the
cyclic group of order $m$ (resp. of order $\infty $) with generator $x$.

\bigskip

\noindent \textbf{Case 1.} $G=SU(n)$. For each $1\leq k\leq n-1$ we have by
(4.8) the exact sequence

\begin{quote}
$0\rightarrow \mathbb{Z}\{\rho _{k}^{\prime }\}\overset{a_{k}}{\rightarrow }%
\mathbb{Z}\{\rho _{k}\}\rightarrow \mathbb{Z}_{b_{n,k}}\{\omega ^{k}\}%
\overset{\omega }{\rightarrow }\omega \cdot \mathbb{Z}_{b_{n,k}}\{\omega
^{k}\}\rightarrow 0$,
\end{quote}

\noindent where $A^{2k}(PG)=\mathbb{Z}_{b_{n,k}}\{\omega ^{k}\}$ by (2.6),
and where

\begin{quote}
$E_{3}^{2k,1}(PG)/\left\langle \omega \right\rangle =\mathbb{Z}\{\rho
_{k}^{\prime }\}$, $E_{3}^{2k,1}(\widetilde{G})=\mathbb{Z}\{\rho _{k}\}$ by
(4.9).
\end{quote}

\noindent As the order of the power $\omega ^{r}$ is precisely $b_{n,r}$, we
get by the exactness that $\ker \omega =\mathbb{Z}_{b_{n,k}/b_{n,k+1}}\{%
\omega \}$, showing $a_{k}=\frac{b_{n,k}}{b_{n,k+1}}$, as that stated in
Theorem D.

\bigskip

\noindent \textbf{Case 2.} $G=Sp(n)$ with $n=2^{r}(2b+1)$. Taking in (4.8)
that $k=2^{r+1}-1$ we get the exact sequence

\begin{quote}
$0\rightarrow \mathbb{Z}\{\rho _{2^{r}}^{\prime }\}\overset{a_{2^{r}}}{%
\rightarrow }\mathbb{Z}\{\rho _{2^{r}}\}\overset{\theta }{\rightarrow }%
\mathbb{Z}_{2}\{\omega ^{2^{r+1}-1}\}\overset{\omega }{\rightarrow }\omega
\cdot \mathbb{Z}_{2}\{\omega ^{2^{r+1}-1}\}=0$,
\end{quote}

\noindent where $\mathcal{A}^{2k}(PG)=\mathbb{Z}_{2}\{\omega ^{k}\}$ and $%
\omega \cdot \mathbb{Z}_{2}\{\omega ^{k}\}=0$ by (2.7), and where

\begin{quote}
$E_{3}^{2k,1}(\widetilde{G})=\mathbb{Z}\{\rho _{2^{r}}\}$, $%
E_{3}^{2k,1}(PG)/\left\langle \omega \right\rangle =\mathbb{Z}\{\rho
_{2^{r}}^{\prime }\}$ by (4.9).
\end{quote}

\noindent This shows that $a_{2^{r}}=2$. Moreover, with $\left\vert \mathcal{%
Z}(G)\right\vert =2$ the relation (4.10) forces $a_{i}=1$ for $i\neq 2^{r}$,
verifying Theorem D for $G=Sp(n)$.

\bigskip

\noindent \textbf{Case 3.} $G=E_{7}$. Taking $k=1$ in (4.8) we obtain the
short exact sequence

\begin{quote}
$0\rightarrow \mathbb{Z}\{\rho _{1}^{\prime }\}\overset{a_{1}}{\rightarrow }%
\mathbb{Z}\{\rho _{1}\}\overset{\theta }{\rightarrow }\mathbb{Z}_{2}\{\omega
\}\overset{\omega }{\rightarrow }\omega \cdot \mathbb{Z}_{2}\{\omega \}=0$
\end{quote}

\noindent where $\mathcal{A}^{2}(PG)=\ker \omega =\mathbb{Z}_{2}\{\omega \}$
by (2.11), and where

\begin{quote}
$E_{3}^{2,1}(\widetilde{G})=\mathbb{Z}\{\rho _{1}\}$, $E_{3}^{2,1}(PG)/\left%
\langle \omega \right\rangle =\mathbb{Z}\{\rho _{1}^{\prime }\}$ by (4.9).
\end{quote}

\noindent It shows that $a_{1}=2$. Consequently, $a_{2}=\cdots =a_{7}=1$ by
(4.10).

\bigskip

\noindent \textbf{Case 4.} $G=E_{6}$. Taking $k=8$ in (4.8) we get the short
exact sequence

\begin{quote}
$0\rightarrow E_{3}^{16,1}(PG)/\left\langle \omega \right\rangle \overset{%
C^{\ast }}{\rightarrow }E_{3}^{16,1}(\widetilde{G})\overset{\theta }{%
\rightarrow }\mathbb{Z}_{3}\{\omega ^{8}\}\rightarrow 0$
\end{quote}

\noindent where $\mathcal{A}^{16}(PG)=\mathbb{Z}_{3}\{\omega ^{8}\}\oplus 
\mathbb{Z}_{3}\{x_{4}^{2}\}$ with $\ker \omega =\mathbb{Z}_{3}\{\omega
^{8}\} $ by (2.10), and where

\begin{quote}
$E_{3}^{16,1}(\widetilde{G})=\mathbb{Z}\{\rho _{5}\}\oplus x_{3}\cdot 
\mathbb{Z}\{\rho _{3}\}\oplus x_{4}\cdot \mathbb{Z}\{\rho _{2}\}$,

$E_{3}^{16,1}(PG)/\left\langle \omega \right\rangle =\mathbb{Z}\{\rho
_{5}^{\prime }\}\oplus x_{3}\cdot \mathbb{Z}\{\rho _{3}^{\prime }\}\oplus
x_{4}\cdot \mathbb{Z}\{\rho _{2}^{\prime }\}$ by (4.9).
\end{quote}

\noindent Moreover, by the relation i) of Theorem 4.2, together $\func{Im}%
\theta =\mathbb{Z}_{3}\{\omega ^{8}\}$,

\begin{quote}
$\theta (x_{3}\cdot \rho _{3})=x_{3}\cdot \theta (\rho _{3})=0$; $\theta
(x_{4}\cdot \rho _{2})=x_{4}\cdot \theta (\rho _{2})=0$.
\end{quote}

\noindent These imply, by the exactness, that $\theta (\rho _{5})=\omega
^{8} $. Since $\theta (\rho _{5})$ is of order $3$ there must be $C^{\ast
}(\rho _{5}^{\prime })=3\rho _{5}$, showing $a_{5}=3$. As result $a_{i}=1$
for $i\neq 5$ by (4.10).

\bigskip

\noindent \textbf{Case 5. }$G=Spin(2n+1)$, $2^{s}\leq n<2^{s+1}$. By the
formula (2.8) of $\mathcal{A}^{\ast }(PG)$

\begin{quote}
$\ker \{A^{2k}(PG)\overset{\omega }{\rightarrow }\omega \cdot
A^{2k}(PG)\}=\left\{ 
\begin{tabular}{l}
$0$ if $0\leq k<2^{s+1}-1$; \\ 
$\mathbb{Z}_{2}\{\omega ^{2^{s+1}-1}\}$ if $k=2^{s+1}-1$,%
\end{tabular}%
\right. $
\end{quote}

\noindent where $\omega =x_{1}$. By (4.8) we have the exact sequences

\begin{enumerate}
\item[(4.12)] $0\rightarrow E_{3}^{2k,1}(PG)/\left\langle \omega
\right\rangle \overset{C^{\ast }}{\rightarrow }E_{3}^{2k,1}(\widetilde{G})%
\overset{\theta }{\rightarrow }0$ for $1\leq k<2^{s+1}-1$;

\item[(4.13)] $0\rightarrow E_{3}^{2(2^{s+1}-1),1}(PG)/\left\langle \omega
\right\rangle \overset{C^{\ast }}{\rightarrow }E_{3}^{2(2^{s+1}-1),1}(%
\widetilde{G})\overset{\theta }{\rightarrow }\mathbb{Z}_{2}\{\omega
^{2^{s+1}-1}\}\rightarrow 0.$
\end{enumerate}

\noindent In particular, by (4.12)

\begin{enumerate}
\item[(4.14)] $\theta (\rho _{i})=0$, $1<i<2^{s}$.
\end{enumerate}

\noindent Moreover, by (4.9) the groups $E_{3}^{2(2^{s+1}-1),1}(\widetilde{G}%
)$ and $E_{3}^{2(2^{s+1}-1),1}(PG)/\left\langle \omega \right\rangle $ in
(4.13) are spanned respectively by the following elements

\begin{quote}
$\rho _{2^{s}}$, $a_{i}\cdot \rho _{i},$ $1\leq i<2^{s}$, $a_{i}\in \mathcal{%
A}^{4(2^{s}-i)}(G)$ and

$\rho _{2^{s}}^{\prime },$ $a_{i}\cdot \rho _{i}^{\prime },$ $1\leq i<2^{s}$%
, $a_{i}\in \mathcal{A}^{4(2^{s}-i)}(G)$,
\end{quote}

\noindent with respect to them, by i) of Lemma 4.2 and by (4.14),

\begin{quote}
$\theta (a_{i}\cdot \rho _{i})=a_{i}\cdot \theta (\rho _{i})=0$, $1\leq
i<2^{s}$,
\end{quote}

\noindent By the exactness of the sequence (4.13) we get

\begin{quote}
$\theta (\rho _{2^{s}})=\omega ^{2^{s+1}-1}$, hence $C^{\ast }(\rho
_{2^{s}}^{\prime })=2\rho _{2^{s}}$,
\end{quote}

\noindent showing $a_{2^{s}}=2$. With $\mathcal{Z}(G)=\mathbb{Z}_{2}$ we get 
$a_{i}=1$ for $i\neq 2^{s}$ by (4.10).

\bigskip

\noindent \textbf{Case 6. }$G=Spin(2n)$, $2^{s}<n=2b+1\leq 2^{s+1}$.
According to (2.9) we have

\begin{quote}
$\mathcal{A}^{\ast }(PG)=\frac{\mathbb{Z[}\omega ,x_{3},x_{5},\cdots ,x_{2%
\left[ \frac{n}{2}\right] -1}]}{\left\langle 4\omega ,2\omega ^{2},\omega
^{2^{s+1}},2x_{2i-1},\text{ }x_{2i-1}^{k_{i}}\text{; }2\leq i\leq \left[ 
\frac{n}{2}\right] \right\rangle }$,
\end{quote}

\noindent where $\omega =x_{1}$. It implies that

\begin{quote}
$\ker \{A^{2k}(PG)\overset{\omega }{\rightarrow }\omega \cdot
A^{2k}(PG)\}=\left\{ 
\begin{tabular}{l}
$\mathbb{Z}_{2}\{2\omega \}$ if $k=1$; \\ 
$\mathbb{Z}_{2}\{\omega ^{2^{s+1}-1}\}$ if $k=2^{s+1}-1$.%
\end{tabular}%
\right. $
\end{quote}

\noindent By (4.8) we have the exact sequences

\begin{quote}
$0\rightarrow \mathbb{Z\{\rho }_{1}^{\prime }\mathbb{\}}\overset{C^{\ast }}{%
\rightarrow }\mathbb{Z\{\rho }_{1}\mathbb{\}}\overset{\theta }{\rightarrow }%
\mathbb{Z}_{2}\{2\omega \}\rightarrow 0$ and

$0\rightarrow E_{3}^{2(2^{s+1}-1),1}(PG)/\left\langle \omega \right\rangle 
\overset{C^{\ast }}{\rightarrow }E_{3}^{2(2^{s+1}-1),1}(\widetilde{G})%
\overset{\theta }{\rightarrow }\mathbb{Z}_{2}\{\omega
^{2^{s+1}-1}\}\rightarrow 0$
\end{quote}

\noindent implying respectively that $a_{1}=2$ and $a_{2^{s}}=2$ (by an
argument similar to the previous case). With $\mathcal{Z}(G)=\mathbb{Z}_{4}$
we get $a_{i}=1$ for $i\neq 1$ or $2^{s}$ by (4.10).

\bigskip

\noindent \textbf{Case 7.} $G=Spin(2n)$, $2^{s}<n=2^{t}(2b+1)\leq
2^{s+1},t\geq 1$. Decompose the covering $c:G\rightarrow PG$ into the
composition of two cyclic coverings of order $2$

\begin{quote}
$c=c_{2}\circ c_{1}:Spin(2n)\overset{c_{1}}{\rightarrow }SO(2n)\overset{c_{2}%
}{\rightarrow }PSpin(2n)$,
\end{quote}

\noindent where, in addition to formulae (2.9) of $\mathcal{A}^{\ast }(G)$
for $G=Spin(2n)$, $PSpin(2n)$,

\begin{quote}
$\mathcal{A}^{\ast }(SO(2n))=\frac{\mathbb{Z[}x_{1},x_{3},x_{5},\cdots ,x_{2%
\left[ \frac{n}{2}\right] -1}]}{\left\langle 2x_{2i-1},\text{ }%
x_{2i-1}^{h_{i}};\text{ }1\leq i\leq \left[ \frac{n}{2}\right] \right\rangle 
}$, $\deg x_{k}=2k,h_{i}=2^{[\ln \frac{n-1}{i}]+1}$.
\end{quote}

\noindent by Marlin \cite{M}. The same calculation as that in Case 5 shows
that

\begin{quote}
$\mathcal{D}(Spin(2n),SO(2n))=\{1,\cdots ,1,2_{(2^{s})},1,\cdots ,1\}$,

$\mathcal{D}(SO(2n),PSpin(2n))=\{1,\cdots ,1,2_{(2^{t-1})},1,\cdots ,1\}$.
\end{quote}

\noindent For $G=Spin(2n)$ Theorem D is verified by the multiplicative
property of multi--degree with respect to the composition of coverings.$%
\square $

\bigskip

For a cyclic covering $c:G\rightarrow G^{\prime }$ on simple Lie groups let $%
J(\omega )\subseteq H^{\ast }(G^{\prime })$ be the subring generated by the
Euler class $\omega $ of the central extension $\widetilde{G}\rightarrow
G^{\prime }$, and let $b_{r}$ be the order of the power $\omega ^{r}$ in $%
J(\omega )$. Define \textsl{the characteristic of the ring }$J(\omega )$ to
be the sequence $Ch(J(\omega )):=\left\{ b_{1},b_{2},\cdots \right\} $ of
integers. It is clear that

\begin{quote}
i) $b_{r+1}$ divides $b_{r}$ for $r\geq 1$; and

ii) $J(\omega )=J(PG)$ for $G\neq Spin(2n)$ with $n\equiv 1\func{mod}2$
\end{quote}

\noindent by Theorem C. For examples, we get from ii) that

\begin{quote}
$Ch(J(\omega ))=\left\{ 3,\cdots ,3,1\right\} $ for the three sheets
covering $E_{6}\rightarrow PE_{6}$

$Ch(J(\omega ))=\left\{ 2,1,\cdots ,1\right\} $ for the two sheets covering $%
E_{7}\rightarrow PE_{7}$.
\end{quote}

\noindent Since the proof of Theorem D goes through (essentially) all the
possible cyclic coverings between simple Lie groups, it implies the
following result that is independent of the types of simple Lie groups.

\bigskip

\noindent \textbf{Theorem 4.4. }\textsl{For any cyclic covering }$%
c:G\rightarrow G^{\prime }$\textsl{\ between simple Lie groups the two
invariants }$\mathcal{D}(G,G^{\prime })$\textsl{\ and }$Ch(J(\omega ))$%
\textsl{\ are related by the equalities }

\begin{quote}
$a_{k}=b_{r_{k}}/b_{r_{k}+1}$\textsl{, }$1\leq k\leq n$, \textsl{where }$%
I_{G}=\{r_{1},\cdots ,r_{n}\}$\textsl{.}$\square $
\end{quote}

\bigskip

\noindent \textbf{Remarks 4.5. }The passage from a cyclic covering $%
c:G\rightarrow G^{\prime }$ to its central extension $\widetilde{G}%
\rightarrow G^{\prime }$ is an useful construction in geometry. As examples,
for the cyclic covering $SU(n)\rightarrow PSU(n)$ of order $n$ we have $%
\widetilde{G}=U(n)$, the unitary group of rank $n$; for the cyclic covering $%
Spin(n)\rightarrow SO(n)$ of order $2$ we get the spin$^{c}$ group $%
\widetilde{G}=Spin^{c}(n)$.

In the context of Schubert calculus the integral cohomologies of all simply
connected simple Lie groups $G$ have been constructed in \cite{DZ2}. In \cite%
{D0,D} the exact sequence (4.5) has been applied to extend the construction
to the integral cohomology of the adjoint Lie groups $PG$.$\square $

\end{document}